\documentclass{amsart} %[envcountsame]{llncs}
\pdfoutput=1

\usepackage{epsfig}
\usepackage{amsfonts}
\usepackage{amsmath}
\usepackage{enumerate}
\usepackage{amsthm}
\usepackage[tight]{subfigure}
\usepackage{footmisc}
\usepackage{url}

\newtheorem{theorem}{Theorem}%[section]
\newtheorem{definition}[theorem]{Definition}

\newtheorem{lemma}[theorem]{Lemma}
\newtheorem{prob}[theorem]{Problem}
\newtheorem{corollary}[theorem]{Corollary}
\theoremstyle{remark}
\newtheorem{example}[theorem]{Example}

\newcounter{c-save-lemma}
\newcounter{c-save-main}

\providecommand{\qed}{\hfill$\Box$\vspace{0.2cm}}
\def \T {\mathcal{T}}
\def \proof {\noindent{\bf Proof}\quad}
\def \ident {\leftrightarrow}
\providecommand{\face}[1]{\textit{#1}}
\providecommand{\tw}{\mathrm{tw}}
\providecommand{\vfpg}[2]{\textsc{#1-admissibility(#2)}}

\title{Fixed parameter tractable algorithms in combinatorial topology}

\author{Benjamin A. Burton  \and
  William Pettersson
}
\thanks{The first author is supported by the Australian Research
Council (DP1094516, DP110101104).}

\begin{document}
\maketitle

\begin{abstract}
To enumerate 3-manifold triangulations with a given property, one typically begins
with a set of potential face pairing graphs (also known as dual 1-skeletons),
and then attempts to flesh each graph out into full triangulations using
an exponential-time enumeration.
However, asymptotically most graphs do not result in \emph{any} 3-manifold
triangulation, which leads to significant ``wasted time'' in
topological enumeration algorithms. Here we give a new algorithm to determine
whether a given face pairing graph supports any 3-manifold triangulation,
and show this to be fixed parameter tractable in the treewidth of the graph.

We extend this result to a ``meta-theorem'' by defining a broad class of
properties of triangulations,
each with a corresponding fixed parameter tractable existence algorithm.
We
explicitly implement this algorithm in the most generic setting, and we
identify heuristics that in practice are seen to
mitigate the large constants that so often occur in
parameterised complexity, highlighting the practicality of our techniques.
%Lastly, we show how these results can also be applied directly to census
%algorithms to provide significant real-world speed improvements.
\end{abstract}

\section{Introduction}
In combinatorial topology, a triangulated 3-manifold involves
abstract tetrahedra whose faces are identified or ``glued'' in pairs.
Many research questions involve
looking for a triangulated manifold
which fits certain requirements, or is pathologically bad for certain
algorithms, or breaks some conjecture.
One invaluable tool for such tasks is an exhaustive \emph{census} of
triangulated 3-manifolds.

The first of these was the
census of cusped hyperbolic 3-manifold triangulations on $\leq 5$ tetrahedra by
Hildebrand and Weeks \cite{Hildebrand1989} in 1989, later extended to
$\leq 9$ tetrahedra \cite{Burton2014Cusped,Callahan1999,Thistlethwaite10-cusped8}.
Another much-used example is the
census
of closed orientable prime minimal triangulations of $\leq 6$ tetrahedra by
Matveev \cite{Matveev1998}, later extended to $\leq 12$
tetrahedra \cite{Martelli2001,Matveev2007AlgorithmicTopology}.
%More recent work
%focuses on classes
%of 3-manifold triangulations with desirable properties, such as 1-vertex
%triangulations \cite{Bacher2001,Burton2011PachnerGraphConference}, angle structures
%\cite{Burton2012TautAngle,Futer2011,Hodgson2011} and discrete Morse theory
%\cite{Burton2013Morse,Joswig2004,Joswig2006,Lewiner2003}.

%These tasks are often approached by testing all manifolds on $n$ tetrahedra
%before testing any manifolds on $n+1$ tetrahedra.
In all of these prior works, the authors enumerate all triangulated manifolds on
$n$ tetrahedra by first enumerating all 4-regular multigraphs on $n$ nodes
(very fast),
and then for each graph $G$ essentially modelling every possible triangulation with
$G$ as its dual graph (very slow).
%In practice this involves finding all possible 4-regular multigraphs, and then
%for each graph $G$ in such a set,
%building up every possible triangulation on $G$. 
%This method allows the use of automorphisms of the graph $G$ to reduce the
%running time of the enumeration process.
If any such triangulation built from $G$ is the triangulation of a
3-manifold, we say that $G$ is {\em admissible}.
% , and that $G$ \emph{admits} a triangulation.
If $G$ admits a 3-manifold triangulation with some particular property $p$,
we say that $G$ is {\em p-admissible}.

%\footnote{Other properties of 3-manifolds are often added to this
%definition, depending on the exact problem at hand.}, we say that $G$ is {\em admissible}.
%Often the number of
%admissible graphs is far fewer than the number of 4-regular multigraphs on a
%given number of vertices.

Using state-of-the-art public software \cite{Regina},
%This is not a purely theoretical barrier;
generating such a census on 12 tetrahedra takes 1967 CPU-days, of which
over 1588 CPU-days is spent analysing non-admissible graphs.
Indeed, for a typical census on $\leq 10$ tetrahedra,
%of closed orientable prime minimal triangulations
less than $1\%$ of $4$-regular graphs are admissible \cite{Burton2007}.
Moreover,
Dunfield and Thurston \cite{Dunfield2006} show that the probability
of a random 4-regular graph being admissible tends toward zero as
the size of
the graph increases.
Clearly an efficient method of determining whether a given graph is admissible
could have significant effect on the (often enormous) running time
required to generating such a census.

We use parameterized complexity \cite{Downey1999} to address this issue.
A problem is {\em fixed parameter tractable}
if, when some parameter of the input is fixed, the problem can be
solved in polynomial time in the input size.
In Theorem~\ref{theorem:tw-algo-plain} we show that
to test whether a graph $G$ is admissible is
fixed parameter tractable, where the parameter is the
treewidth of $G$. Specifically, if the
treewidth is fixed at $\leq k$ and $G$ has size $n$, we can
determine whether $G$ is admissible in $O(n\cdot f(k))$ time.

Courcelle showed \cite{Courcelle1990,Courcelle2001} that for graphs of bounded
treewidth, an entire class of problems have fixed parameter tractable algorithms.
However, employing this result for our problem of testing admissibility
looks to be highly non-trivial.
In particular,
it is not clear how the topological constraints of our problem can be
expressed in monadic second-order logic, as Courcelle's theorem requires.
%Courcelle's
%Theorem also imposes significant constants in the parameterized algorithms
%which could cause practical problems.
Even if Courcelle's theorem could be used, our results here provide
significantly better constants than a direct application of
Courcelle's theorem would.

%The use of monadic second-order logic also causes issues. TODO.
%More recent work \cite{Burton2013Morse,Burton2012TautAngle} has also shown that
%for various extra properties, determining whether a graph is admissible is also
%fixed parameter tractable.

Following the example of Courcelle's theorem, however,
we generalise our result to a larger class of problems
(Theorem \ref{theorem:tw-algo}).
%a meta-theorem which provides a similar result
%for closed 3-manifold triangulations on face pairing
%graphs of bounded treewidth.%, which extends these existing results.
Specifically, we introduce the concept of a
{\em simple property}, and give a fixed parameter tractable algorithm
which, for any simple property $p$, determines whether a graph admits a
triangulated 3-manifold with property $p$
(again the parameter is treewidth).

We show that these results are practical through an explicit
implementation, and identify some simple heuristics
which improve the running time
and memory requirements.
To finish the paper, we identify a clear potential for
how these ideas can be extended to the more difficult enumeration
problem, in those cases
where a graph \emph{is} admissible and a
complete list of triangulations is required.

% of our implementation, and can also be applied
% to similar algorithms on tree decompositions.
% this algorithm, and also show how the results can be utilised by existing
% census algorithms to achieve running time improvements of almost 400\%.
%These results are easily extended to other properties, and we also demonstrate
%how to extend the algorithm.

Parameterised complexity is very new to the field of
3-manifold topology
\cite{burton13-morse,Burton2012TautAngle}, and this paper marks the
first exploration of parameterised complexity in 3-manifold enumeration
problems.  Given that 3-manifold algorithms are often extremely slow and
complex, our work here highlights a growing potential for parameterised
complexity to offer practical alternative algorithms in this field.

\section{Background}

To avoid ambiguity with the words ``vertex'' and ``edge'',
we use the terms \emph{node} and \emph{arc} instead for graphs,
and \emph{vertex} and \emph{edge} in the context of triangulations.

% \subsection{Graphs, tree decompositions and treewidth}

Many NP-hard problems on graphs are fixed parameter tractable in the
\emph{treewidth} of the graph
(e.g.,
\cite{Arnborg1985,Arnborg1991,Bodlaender1996,Bodlaender1996Efficient,Courcelle2001}).
Introduced by Robertson and Seymour \cite{Robertson1986},
the treewidth measures precisely how ``tree-like'' a graph is:

\begin{definition}[Tree decomposition and treewidth]
  \label{definition:treedecomp}
  Given a graph $G$, a tree decomposition of $G$ is a tree $H$ with the
  following additional properties:
  \begin{itemize}
    \item Each node of $H$, also called a \emph{bag}, is associated with a set of nodes of $G$;
    \item For every arc a of $G$, some bag of $H$ contains
    both endpoints of $a$; 
    \item For any node $v$ of $G$, the subforest in $H$ of
      bags containing $v$ is connected.
  \end{itemize}
If the largest bag of $H$ contains $k$ nodes of $G$, we say that the
tree decomposition has \emph{width} $k+1$.
%Clearly a graph can have many tree decompositions, so we say that
The \emph{treewidth} of $G$, denoted
$\tw(G)$, is the minimum width of any tree decomposition of $G$.
\end{definition}
%
%\begin{figure}
%  \centering
%  \subfloat[][]{
%    \includegraphics{TreeDecomp-G}
%  }
%  \subfloat[][]{
%    \includegraphics{TreeDecomp-H}
%  }
%
%  \caption{A graph $G$ on the left, and a tree decomposition of $G$ on the
%  right}
%  \label{fig:TreeDecomp}
%\end{figure}

Bodlaender \cite{Bodlaender1996} gave
a linear
time algorithm for determining if a graph has treewidth $\leq k$ for fixed $k$,
and for finding such a tree decomposition,
and Kloks \cite{Kloks1994Treewidth}
demonstrated algorithms for finding ``nice'' tree decompositions.

% \subsection{3-manifolds: Notation and definitions}
\label{sec:partial-3-manifolds}
A closed 3-manifold is essentially a topological space in which every point has
some small neighbourhood homeomorphic to $\mathbb{R}^3$.
%Combinatorially it can be described \cite{Moise1952} as a triangulation; a set of tetrahedra,
%and identifications between faces of said tetrahedra.
We first define {\em
general triangulations}, and then give conditions under which they
represent 3-manifolds.
\begin{definition}[General triangulation]
A general triangulation
%$(\mathcal{D},\mathcal{P})$
is a set of abstract tetrahedra
%$\mathcal{D} =
$\{\Delta_1,\Delta_2,\ldots,\Delta_n\}$ and a set of face
identifications or ``gluings''
%$\mathcal{P} =
$\{\pi_1,\pi_2,\ldots,\pi_m\}$, such that each
$\pi_i$ is an affine identification between two distinct faces of tetrahedra,
and each face is a part of at most one such identification.
\end{definition}

Note that this is more general than a simplicial complex
(e.g., we allow an identification between two distinct faces of the
same tetrahedron), and it
need not represent a 3-manifold.
Any face which is not identified to
another face is called a \emph{boundary face} of the triangulation. If a
triangulation has no such boundary faces, we say it is \emph{closed}. We also
note that there are six ways to identify two faces, given by the six symmetries
of a regular triangle.

We can partially represent a triangulation by its face pairing graph,
% also known as its dual 1-skeleton. This graph
which describes \emph{which} faces are
identified together, but not \emph{how} they are identified.
\begin{definition}[Face pairing graph]
  The face pairing graph of a triangulation $\T$ is
  the multigraph $\Gamma(\T)$ constructed as follows.  Start with an empty graph $G$, and insert
  one node for every tetrahedron in $\T$.  For every face identification
  between two tetrahedra $T_i$ and $T_j$, insert the arc $\{i,j\}$ into the
  graph $G$.
  %We denote the face pairing graph of a triangulation $\T$ by $G = \Gamma(\T)$.
\end{definition}

Note that a face pairing graph will have parallel arcs if there are two
distinct face
identifications between $T_i$ and $T_j$, and loops if two faces of the same
tetrahedron are identified together.

Some edges of tetrahedra will be
identified together as a result of these face identifications
(and likewise for vertices).
Some edges may be identified directly via a single face
identification, while others may be identified indirectly through
a series of face identifications.

%\begin{notation}[Edge identification and orientation]
We assign an arbitrary orientation to each edge of each tetrahedron. 
Given two tetrahedron edges $e$ and $e'$ that are identified together via the face
identifications, we write $e \simeq e'$ if %they are identified so that
the orientations agree, and $e \simeq \overline {e'}$ if the orientations are 
reversed.  In settings where we are not interested in orientation, we write
$e \sim e'$ if the two edges are identified (i.e. one of $e \simeq e'$ or $e
\simeq \overline {e'}$ holds).  
%\end{notation}

This leads to the natural notation $[e] = \{e' : e \sim e'\}$ as an
equivalence class of identified edges (ignoring orientation).
We refer to $[e]$ as an \emph{edge of the triangulation}.
%, and $e$ will refer to an individual edge of a tetrahedron.
Likewise, we use the notation $v \sim v'$ for vertices of tetrahedra
that are identified together via the face identifications,
and we call an equivalence class $[v]$ of
identified vertices a \emph{vertex of the triangulation}.
%together in said triangulation, while $v$ and $v'$ are individual vertices of
%tetrahedra.

%We will also discuss the various boundary components of triangulations.
A {\em boundary edge\,/\,vertex} of a triangulation is an edge\,/\,vertex
of the triangulation
whose equivalence class contains some edge\,/\,vertex of a boundary face.
%Unless otherwise indicated, these terms refer to edges (respectively
%vertices) of the triangulation, not of individual tetrahedra.
%That is, a boundary edge is the equivalence class of edges

%\begin{definition}
%  \textsc{link graph} The link graph of a vertex $[v]$ in a triangulation $\T$ is
%  the graph constructed such that an edge in $\T$ is a node in the link graph
%  of $[v]$ if and only if the edge is incident upon some $v' \in [v]$. Two
%  nodes are adjacent in the link graph if they share a common face which is
%  also incident with some $v' \in [v]$.
%\end{definition}

%Note that if we take a link graph of a vertex $[v]$, and insert a new node and
%ensure this node is adjacent to every other node the resulting graph is
%actually a triangulation of the vertex link.

%The main results in this paper involve finding {\em 3-manifold triangulations}.
%Several definitions will now follow, which will lead to the definition of a
%3-manifold triangulation.

\begin{figure}[t]
  \centering
  %\subfigure[]{\includegraphics[scale=0.5]{gt}\label{fig:gt}}
  %\qquad
  %\qquad
  %\subfigure[]{\includegraphics[scale=0.5]{link}\label{fig:link}}
  \includegraphics[scale=0.5,angle=270]{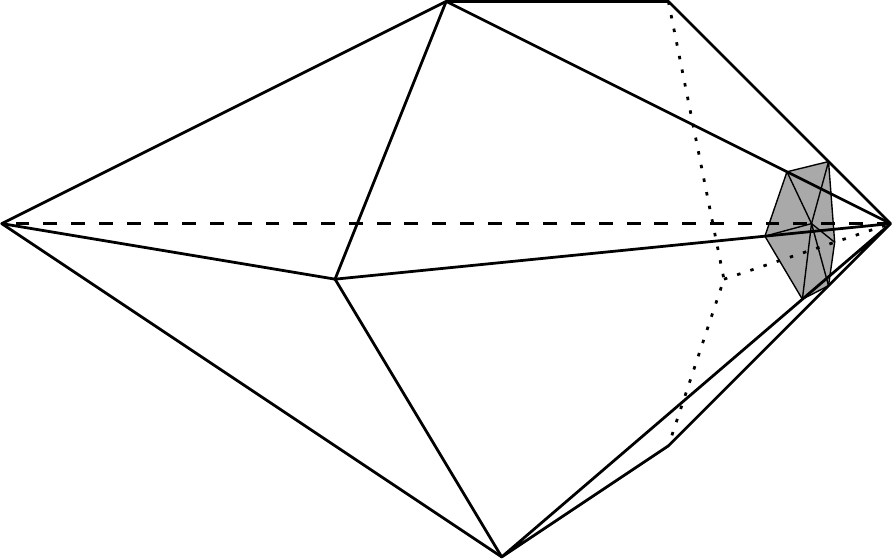}
  \caption{A triangulation of a 3-ball with 6
  tetrahedra meeting along an internal edge.}
  \label{fig:link}
  %Figure \ref{fig:link} shows the link of the top vertex shaded.}
\end{figure}

%\begin{definition}[Vertex link]
  The \emph{link} of a vertex $[v]$ is the
  (2-dimensional) frontier of a small regular neighbourhood of $[v]$.
%\end{definition}
Figure \ref{fig:link} shows the link of the top vertex shaded in
grey; in this figure, the link is homeomorphic (topologically
equivalent) to a disc. The link is
a 2-dimensional triangulation (in the example it has six triangles),
and we use the term {\em arc} to denote an
edge in this triangulation.
In this paper, whether ``arc'' refers to a graph or a
vertex link is always clear from context.

%Below we define two different types of triangulations, both of which place
%restrictions of the links of vertices in the triangulation.

\begin{definition}[Closed 3-manifold triangulation]
A closed 3-manifold triangulation $\T$ is a general triangulation for which
(i) $\T$ is connected; 
(ii) for any vertex $v$ in $\T$, the link of $v$ is homeomorphic to a
2-sphere; and
(ii) no edge $e$ in $\T$ is identified with itself in reverse (i.e. $e \not\simeq \overline e$).
\end{definition}

These properties are necessary and sufficient for the underlying topological
space to be a 3-manifold.
We say that a graph $G$
%{\em admits} a triangulation $\T$ if $G$ is the
%face pairing graph of $\T$. We say that $G$
is {\em admissible} if it is the face pairing graph for any closed
3-manifold triangulation $\T$.

\begin{definition}[Partial-3-manifold triangulation]\label{definition:partial-3mfld}
A partial-3-manifold triangulation $\T$ is a general triangulation for which
(i) for any vertex $v$ in $\T$, the link of $v$ is homeomorphic to a
  2-sphere with zero or more punctures; and
(ii) no edge $e$ in $\T$ is identified with itself in reverse (i.e. $e \not\simeq \overline e$).
\end{definition}

These are in essence ``partially constructed'' 3-manifold triangulations;
the algorithms of Section \ref{sec:algorithm} build these up
into full 3-manifold triangulations.
Note that the underlying space of $\T$ might not even be a 3-manifold
with boundary: there may be ``pinched vertices'' whose links have
many punctures.
%In particular, a 3-manifold triangulation only allows links
%homeomorphic to a 2-sphere with at most one puncture, while in a
%partial-3-manifold triangulation links may have any number of punctures.

We can make some simple observations:
%Some simple observations about partial 3-manifold triangulations:
%\begin{itemize}
(i) the boundary vertices of a partial 3-manifold triangulation
are precisely those whose links have at least one puncture;
  %Note that vertices with links homeomorphic to a
%2-sphere cannot lie on any boundary of a given triangulation, and vertex
%links homeomorphic to a 2-sphere with one more punctures must be boundary
%vertices.
  (ii) a connected partial-3-manifold triangulation with no boundary faces is
    a closed 3-manifold triangulation, and vice-versa;
  (iii) a partial-3-manifold triangulation with a face identification removed,
  or an entire tetrahedron removed,
  is still a
    partial-3-manifold triangulation.

\section{Configurations} \label{sec:config}

The algorithms in Section \ref{sec:algorithm} build up 3-manifold triangulations one
tetrahedron at a time. As we add tetrahedra, we must track what
happens on the boundary of the triangulation, but we can forget about the
parts of the triangulation not on the boundary---this is key
to showing fixed parameter tractability.
In this section we define and analyse
edge and vertex configurations of general
triangulations, which encode exactly those details on the boundary
that we must retain.
%We will shortly
%give an example, and later show that we can
%bound the possible number of configurations of a
%partial-3-manifold-triangulation to a function of the number of boundary faces
%of $\T$.
%Theorem \ref{theorem:tw-algo-plain}.
%, which will give the parametric complexity result.

\begin{definition}[Edge configuration]
The edge configuration of a triangulation $\T$ is a set $C_e$ of
triples detailing how the edges of the boundary faces are identified
together.
Each triple is of the form $((f,e),\ (f',e'),\ o)$, where:
$f$ and $f'$ are boundary faces; $e$ and $e'$ are
tetrahedron edges that lie in $f$ and $f'$ respectively;
$e$ and $e'$ are identified in $\T$;
and $o$ is a boolean ``orientation indicator''
that is true if $e \simeq e'$ and false if $e \simeq \overline{e'}$.
\end{definition}

This mostly encodes the 2-dimensional triangulation of the boundary,
though additional information describing ``pinched vertices'' is still
required.
%by indicating how the edges of boundary faces are identified together in pairs.
%partitioned into pairs, and also includes an orientation marker for each such
%pairing.

\begin{example}[2-tetrahedra pinched pyramid] \label{ex:2tet-edge}
In all examples, we use the
notation $t_i\!:\!a$ to denote vertex $a$ of tetrahedron $t_i$, and
$t_i\!:\!\face{abc}$ to
denote face $\face{abc}$ of tetrahedron $t_i$. Face identifications are denoted as
$t_i\!:\!\face{abc} \ident t_j\!:\!\face{def}$, which means that face $\face{abc}$ of
$t_i$ is mapped to face $\face{def}$ of $t_j$ such that $a\ident d$, $b\ident e$ and
$c\ident f$.

Take two tetrahedra $t_0$ and $t_1$, each with vertices labelled $0,1,2,3$,
and apply the face identifications
  $t_0\!:\!012 \ident t_1\!:\!012$ and
  $t_1\!:\!023 \ident t_1\!:\!321$.

The resulting triangulation is a square based pyramid with one pair of
opposing faces identified (see Figure \ref{fig:2tet-tri}). The final space
resembles a hockey puck with a pinch in the centre, as seen in Figure
\ref{fig:puck}. Note that the vertex at top of the pyramid,
which becomes the pinched centre of the puck, has a link
homeomorphic to a 2-sphere with two punctures. Therefore, although this is a
partial 3-manifold triangulation, the underlying space is not a 3-manifold.

\begin{figure}[t]
  \centering
  \subfigure[]{\includegraphics[scale=0.7]{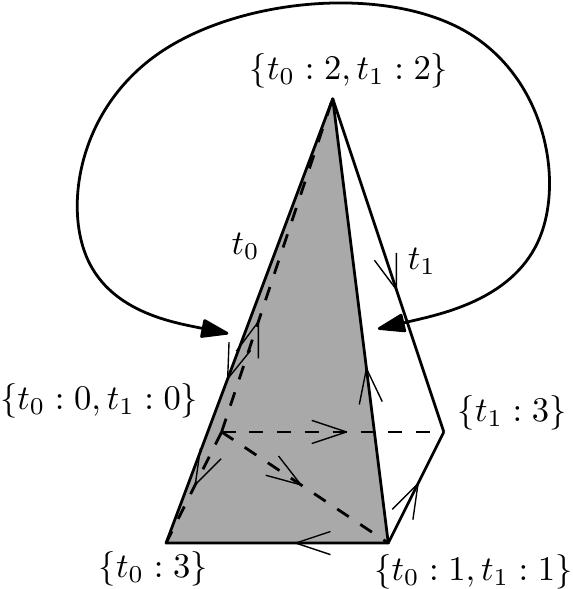}\label{fig:2tet-tri}}
  \qquad
  \subfigure[]{\includegraphics{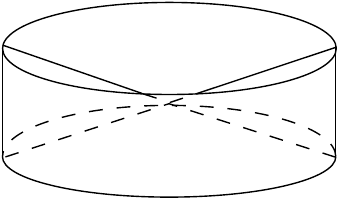}\label{fig:puck}}
  \caption{The triangulation
  from Example \ref{ex:2tet-edge}. The grey shaded tetrahedron is $t_0$.
  Edges are marked with their orientations, and the
  double-ended arrow indicates the identification of two opposing faces of the
  pyramid. The resulting space resembles a
  hockey puck with the centre pinched into a point. This pinch is the vertex
  $\{t_0\!:\!2,t_1\!:\!2\}$.
  }
  \label{fig:2tet}
\end{figure}

The edge configuration of this triangulation is:
\begin{align*}
\{
  &((t_0\!:\!013,03),(t_1\!:\!013,13),f), &(&(t_0\!:\!013,01),(t_1\!:\!013,01),t),& \\
  &((t_0\!:\!013,13),(t_0\!:\!123,13),t), &(&(t_0\!:\!123,12),(t_0\!:\!123,23),f),& \\
  &((t_1\!:\!013,03),(t_1\!:\!023,03),t), &(&(t_1\!:\!023,02),(t_1\!:\!023,23),f)
\}; &
\end{align*}
here $t$ and $f$ represent {\em true} and {\em false} respectively.
%Note how in the first triple in the edge configuration, edge $03$ of face
%$t_0\!:\!013$ is identified to edge $13$ of face $t_1\!:\!013$ in reverse (indicated by
%$f$). If we reverse the orientation on the edge $13$ of face $t_1\!:\!013$ we would
%have written this triple as
%$((t_0\!:\!013,03),(t_1\!:\!013,31,t))$ but we have chosen canonical orientations for
%all edges.
\end{example}
%For each edge $e$ on a boundary face $f$ of $\T$, $E$ contains the tuple
%$(\{(f,e),(f',e')\}, 0)$ if $e \simeq e'$, and $(\{(f,e),(f',e')\}, 1)$ if $e \simeq
%\overline {e'}$.
%Theorem \ref{theorem:ring-bdry} tells us that for each $(f,e)$ pair there exists
%exactly one $(f',e')$ pair that can satisify either of these constraints, and
%exactly one of $e \sim e'$ or $e \sim overline {e'}$ must hold.

\begin{definition}[Vertex configuration]
The vertex configuration $C_v$ of a triangulation $\T$ is a partitioning
of those tetrahedron vertices that
belong to boundary faces, where vertices
$v$ and $v'$ are in the same partition if and only if $v \sim v'$.
\end{definition}

In partial-3-manifold triangulations, vertex links may have
multiple punctures; the vertex configuration then allows us to deduce
which punctures belong to
the same link. In essence, the vertex configuration describes how
the triangulation is ``pinched'' inside the manifold at vertices whose links
have too many punctures.

%\vspace{2mm}
For instance, the vertex configuration of Example \ref{ex:2tet-edge} is given by
\begin{align*}
  \{\{t_0\!:\!0,t_1\!:\!0,t_1\!:\!3\},\quad
  \{t_0\!:\!1,t_0\!:\!3,t_1\!:\!1\},\quad
  \{t_0\!:\!2,t_1\!:\!2\}\}.
\end{align*}
The partition $\{t_0\!:\!2,t_1\!:\!2\}$ represents the pinch at the center of the
``hockey puck''.
%See Appendix~\ref{app:examples} for an additional example (the
%5-tetrahedron cube).
We now give the boundary configuration of a triangulated cube on 5 tetrahedra
as an additional example.

\begin{example}[5-tetrahedra cube]
Take 5 tetrahedra labelled $t_0,t_1,\ldots,t_4$, each with vertices labelled
$0,1,2,3$ and identify the following faces:
\begin{align*}
  t_0\!:\!012 \ident t_4\!:\!012 \\
  t_1\!:\!013 \ident t_4\!:\!013 \\
  t_2\!:\!123 \ident t_4\!:\!123 \\
  t_3\!:\!023 \ident t_4\!:\!023
\end{align*}

The resulting triangulation is a cube, with $t_4$ having zero boundary faces
and each other tetrahedra having three boundary faces. The edge configuration
of this triangulation is:
\begin{align*}
\{
  &((t_0\!:\!013,01),(t_1\!:\!012,01),t), &(&(t_0\!:\!123,12),(t_2\!:\!012,12),t),& \\
  &((t_0\!:\!023,02),(t_3\!:\!012,02),t), &(&(t_1\!:\!123,13),(t_2\!:\!013,13),t),& \\
  &((t_1\!:\!023,03),(t_3\!:\!013,03),t), &(&(t_2\!:\!023,23),(t_3\!:\!123,23),t),& \\
  &((t_0\!:\!013,03),(t_0\!:\!023,03),t), &(&(t_0\!:\!013,13),(t_0\!:\!123,13),t),& \\
  &((t_0\!:\!023,23),(t_0\!:\!123,23),t), &(&(t_1\!:\!023,02),(t_1\!:\!012,02),t),& \\
  &((t_1\!:\!023,23),(t_1\!:\!123,23),t), &(&(t_1\!:\!012,12),(t_1\!:\!123,12),t),& \\
  &((t_2\!:\!023,02),(t_2\!:\!012,02),t), &(&(t_2\!:\!023,03),(t_2\!:\!013,03),t),& \\
  &((t_2\!:\!012,01),(t_2\!:\!013,01),t), &(&(t_3\!:\!012,01),(t_3\!:\!013,01),t),& \\
  &((t_3\!:\!012,12),(t_3\!:\!123,12),t), &(&(t_3\!:\!013,13),(t_3\!:\!123,13),t)
\}. &
\end{align*}

If we examine the first triple in this configuration, it indicates that edge $01$ on
face $013$ of $t_0$  %(in this case through multiple
%identifications)
and edge $01$ on face $012$ of $t_1$ are part of the same edge of the
triangulation\footnote{Remember that each edge of the triangulation is an equivalence
class of edges of tetrahedra.}, even though no face of $t_0$ is
identified to a face of $t_1$.

%The $t$ indicates
%that this identification preserves orientation (that is, vertex $0$ on $t_0$ is
%identified with vertex $0$ on $t_1$ and not vice-versa).

This triangulation has no pinched vertices, but some vertices of tetrahedra are
identified together, and so the vertex configuration is:
\begin{align*}
  \{ \{t_0\!:\!0,t_1\!:\!0,t_3\!:\!0\},\quad
  \{t_0\!:\!1,t_1\!:\!1,t_2\!:\!1\}, \quad\{t_0\!:\!2,t_2\!:\!2,t_3\!:\!2\}, \\
  \{t_0\!:\!3,\}, \quad\{t_1\!:\!3,t_2\!:\!3,t_3\!:\!3\}, \quad\{t_1\!:\!2\},
  \quad\{t_2\!:\!0\}, \quad\{t_3\!:\!1\} \}.
\end{align*}
\end{example}

\begin{definition}[Boundary configuration]
  The boundary configuration $C$ of a triangulation $\T$ is the pair
  $(C_e,C_v)$ where $C_e$ is the edge configuration and $C_v$ is the vertex
  configuration.
\end{definition}

% The boundary configuration is just a combination of the edge and vertex
% configurations. While we have described them separately in this section, we
% will only refer to the combined boundary configurations in section
% \ref{sec:algorithm}.

\begin{lemma}\label{lemma:edge-config-b}
  For $b$ boundary faces,
  there are
  %these boundary
  %faces and their edges can be configured into at most
$\frac{(3b)!}{(3b/2)!}$
%$\frac{(6m)!} {2^{3m-1}\cdot (3m)! }$
possible edge configurations.
\end{lemma}

\proof
Note that $b$ must be even; let $b=2m$.
Each boundary face has three edges, so there are
$6m$ possible pairs $(f,e)$ where $e$ is an edge on a boundary face $f$. Each
such pair must be identified with exactly one other pair, with either
$e
\simeq e'$ or $e \simeq \overline{e'}$, and so the number of possible
edge configurations is
\[2\cdot(6m-1)\cdot 2\cdot (6m-3)\cdot \ldots \cdot 2\cdot 3\cdot 2\cdot 1 =
\frac{(6m)!}{(3m)!} = \frac{(3b)!}{(3b/2)!}. \] \qquad \qed

\begin{lemma}\label{lemma:exact-num-bdry-configs}
For $b$ boundary faces,
the number of possible boundary configurations 
is bounded from above by
\[
  \frac{(3b)!}{(3b/2)!}\cdot \left( \frac{2.376b}{\ln(3b+1)} \right)^{3b}.
\]
%where $B_{3k}$ is the $3k$'th Bell number.

\end{lemma}
\proof
There are
$3b$ tetrahedron vertices on boundary faces, and so the number of possible vertex configurations is
the Bell number $B_{3b}$.  The result now follows from
Lemma~\ref{lemma:edge-config-b} and the following inequality
of Berend \cite{Berend2010}:
\[
  B_{3b} = \frac{1}{e} \sum_{i=0} ^\infty \frac{i^{3b}}{i!} <
  \left( \frac{2.376b}{\ln(3b+1)} \right)^{3b}.
\]\qquad\qed

\begin{corollary}\label{coro:num-bdry-configs}
The number of possible boundary configurations for a triangulation on $n$
tetrahedra with $b$ boundary faces depends on
$b$, but not on $n$.
%bounded by $(9k^2-3k)\cdot B_{3k}$
%where $B_{3k}$ is the $3k$'th Bell number.
\end{corollary}
%\proof
%This follows immediately from the definition of a boundary configuration and
%Lemma \ref{lemma:edge-config-b}. \qed

\begin{figure}[tb]
  \centering
  \includegraphics[scale=0.8]{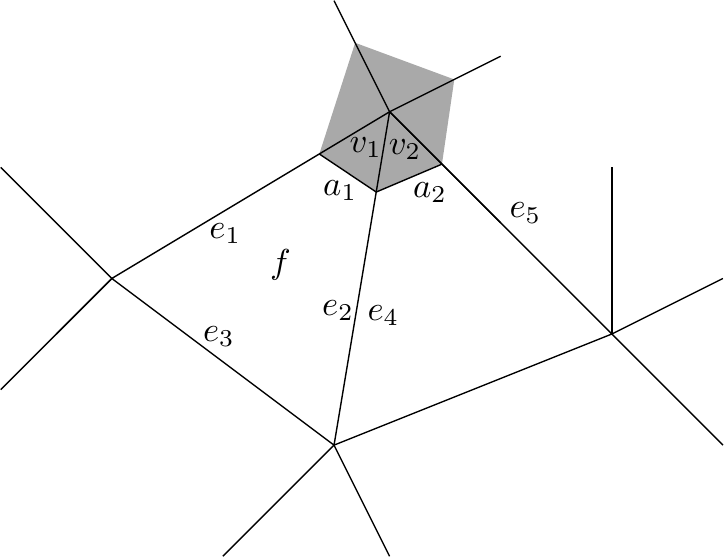}
  \caption{Part of the boundary of a triangulation.
  %One boundary face is marked $f$, with the three edges $e_1$,
  %$e_2$ and $e_3$.
  The link of the top vertex is shaded grey; this link
  does not contain the vertex, but instead cradles the vertex from below.}
  %The arcs $a_1$ and $a_2$ can be reconstructed from the boundary configuration.}
  \label{fig:reconstruct-link}
\end{figure}

The boundary configuration can be used to partially reconstruct the links of
vertices on the boundary of the triangulation.  In particular:
\begin{itemize}
\item The edge configuration allows us to follow the arcs around each puncture
of a vertex link---in Figure~\ref{fig:reconstruct-link} for instance,
we can follow the sequence of arcs $a_1,a_2,\ldots$ that surround the
puncture in the link of the top vertex.
\item The vertex configuration tells us whether two sequences of arcs describe
punctures in the \emph{same} vertex link, versus \emph{different} vertex links.
\end{itemize}

In this way, we can reconstruct all information about punctures in the
vertex links, even though we cannot access the full (2-dimensional)
triangulations of the links themselves.
As the next result shows, this means that the boundary configuration retains all
data required to build up a partial-3-manifold triangulation,
without knowledge of the full triangulation of the underlying space.

%In particular, these links
%must contain some punctures. Recall that the term arc refers to an edge of the
%2-dimensional triangulation of a link. Figure \ref{fig:reconstruct-link} shows
%how the arc $a_1$ can be drawn from $e_1$ to $e_2$, where $f$ is a boundary face,
%$e_1$, $e_2$ and $e_3$ are the edges of $f$ and $v_1$ is a vertex of $f$. The
%arc $a_1$ partially describes the puncture in the link of $[v_1]$ (shaded grey in
%Figure \ref{fig:reconstruct-link}). Note that $e_2$
%must be paired with some other edge in the configuration ($e_4$ in this case),
%and we can also use
%the edge configuration to determine that $v_1$ is identified to $v_2$.
%Therefore the next arc about the puncture must be $a_2$, running from $e_4$ to
%$e_5$ as $v_2$ is incident with both these edges. We can repeat this process to
%reconstruct all arcs about this puncture, and then use the same procedure to
%reconstruct all punctures about vertices on boundary faces.

%We can then use the vertex configuration to determine which vertices of
%boundary faces are part of the same vertex in the underlying space, and
%therefore which punctures are part of the same vertex links.

\begin{lemma}\label{lemma:add-one-id}
Let $\T$ be a partial 3-manifold triangulation
with $b$ boundary faces, and let $\T'$ be formed by
introducing a new identification between two boundary faces of $\T$.
Given the boundary configuration of $\T$ and the new face identification,
we can test whether $\T'$
is also a partial-3-manifold triangulation in $O(b)$ time.
\end{lemma}

%A full proof appears in
%Appendix~\ref{app:add-one-id}. The basic idea is to check whether the
%conditions in Definition~\ref{definition:partial-3mfld} are preserved.
%The edge configuration allows us to easily test for
%edges identified together in reverse,
%and the partial reconstruction of the vertex links
%(as described above) allows us to test whether all
%vertex links are still 2-spheres with zero or more punctures.

%Here we give a full proof of Lemma~\ref{lemma:add-one-id}, which was
%omitted from the main text due to space considerations.
%
%\setcounter{theorem}{\arabic{c-save-lemma}}
%\begin{lemma}
%Let $\T$ be a partial 3-manifold triangulation
%with $b$ boundary faces, and let $\T'$ be formed by
%introducing a new identification between two boundary faces of $\T$.
%Given the boundary configuration of $\T$ and the new face identification,
%we can test whether $\T'$
%is also a partial-3-manifold triangulation in $O(b)$ time.
%\end{lemma}

\proof We are given the boundary
configuration for $\T$. We need to check that the new face identification does
not result in any edges being identified in reverse, and that the links of
all vertices are still spheres with zero or more punctures.

It is easy to see that the face identification will identify at most three pairs
of boundary edges together, and it is routine to check (using the
edge configuration) whether these identifications
will result in any edge identified with itself in reverse in $O(b)$ time.
The rest of this proof therefore only deals with the vertex links.

%As mentioned, each face identification results in three edge identifications,
%and each corresponds to an arc identification in
%the link graph of some vertex of $\T$. In the picture below, we see the three
%edges $ab,bc$ and $ac$ are being identified with $de,ef$ and $df$ respectively. We
%therefore know that arc $\{ab,bc\}$ in the link graph of $v$ is being identified
%with the arc $\{de,ef\}$, and also that in this identification we have $ab \ident
%de$ and $bc \ident ef$.

To determine whether $\T'$ is also a partial-3-manifold triangulation we only
need to determine how these new edge identifications affect the link of each
vertex. Clearly any vertices that are internal to $\T$ must already have links
homeomorphic to a 2-sphere, and cannot be changed. We noted in
Section~\ref{sec:config} that the
cycles of arcs surrounding the punctures on
the links of each boundary vertex can be determined from the
edge configuration---we do not explicitly reconstruct these cycles here,
but we do note that this information is accessible from
the edge configuration.
We note also that each link of a boundary vertex in $\T$
must be homeomorphic to a sphere with
one or more punctures (equivalently, a disc with zero or more punctures).

\begin{figure}[tb]
  \centering
  \includegraphics{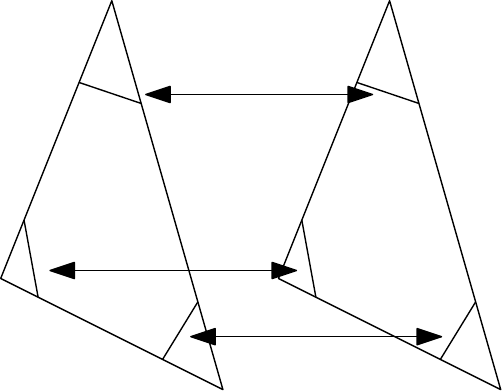}
  \caption{One face identification will result in three pairs of arcs on vertex
  links being identified, as shown by the arrows.}
  \label{fig:arc-identified}
\end{figure}

\begin{figure}
  \centering
  \subfigure[]{\includegraphics[scale=0.8]{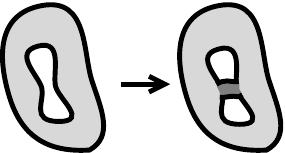}}
  \qquad
  \qquad
  \subfigure[]{\includegraphics[scale=0.8]{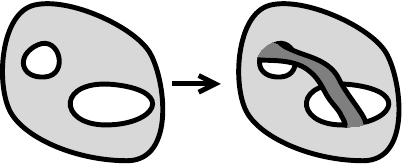}\label{fig:LinkEdgeJoin2}}
  \qquad
  \qquad
  \subfigure[]{\includegraphics[scale=0.8]{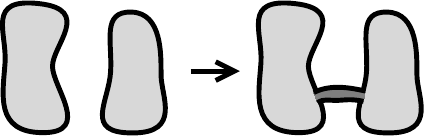}}
  \caption{(a) Type I arc identification; (b) Type II arc identification; and
  (c) Type III arc identification. Light grey indicates the existing link with
  the white space indicating the punctures, and
  the dark grey indicates the new identification in the link.}
  %\caption{From left to right, Type I, Type II and Type III link-edge
  %  identifications. Light grey indicates the existing link, and the dark grey
  %indicates the new identification in the link.}
  \label{fig:LinkEdgeJoins}
\end{figure}

The new face identification will identify three pairs of arcs on the vertex links,
as shown in Figure \ref{fig:arc-identified}. Each of these three arc identifications will take one
of three forms (see Figure \ref{fig:LinkEdgeJoins}):
\begin{enumerate}[Type I:]
  \item The two arcs being identified both bound the same puncture in
  the same vertex link.
  \item The two arcs are part of the same vertex link but bound distinct
    punctures.
  \item The two arcs are part of distinct vertex links.
\end{enumerate}

For a type~I identification,
if the identification preserves the orientability of the
vertex link then the new vertex link will be homeomorphic to a 2-sphere with
zero or more punctures. In particular, if the puncture only
contained two arcs and these are now both identified (in an orientable manner)
then the puncture will be closed off (i.e., the vertex link will be a
2-sphere with one less puncture than before). If this vertex of the triangulation only
had one puncture, then the vertex link will become homeomorphic to a 2-sphere and
the new vertex will be an internal vertex of $\T'$.

If orientability is not preserved in a type~I identification
then we will embed a
M\"obius band in the vertex link, which is never allowed.
Identifications of type II increase the genus of the vertex link,
which is likewise not allowed (see Figure
\ref{fig:LinkEdgeJoin2}), and
identifications of type III simply connect two discs with zero or more
punctures.

In summary: orientable identifications of type I and all identifications of
type III are allowed, whereas non-orientable identifications of type I and all
identifications of type II are not allowed.

Since the triangulations of all vertex links contain $3b$ boundary arcs
in total, we can
identify both the type and orientation of each
identification in $O(b)$ time. Specifically, we
use the edge configuration to determine if the identification is of type~I
(as well as
the orientation of the identification), and we use the vertex configuration to
distinguish between identifications of type~II and type~III. If
any non-orientable type~I identifications or any type II identifications are
found, $\T'$ is not a partial-3-manifold triangulation.

Since we have only
three such identifications of pairs of arcs, we can check all three in $O(b)$
time as well.  Combining this with the $O(b)$ check described earlier
for bad edges, we obtain the required result.  \qed

\section{Algorithms and simple properties}

Recall that the motivating problem for our work was to quickly detect
whether a given graph admits a closed 3-manifold triangulation.
To this end we show:

\begin{theorem}\label{theorem:tw-algo-plain}
%\textsc{valid face pairing graph}
Given a connected 4-regular multigraph $G$, the problem of determining
whether there exists a
closed 3-manifold triangulation $\T$ such that $\Gamma(\T) = G$
is fixed parameter tractable in the treewidth of $G$.
\end{theorem}

This is a special case of our more general Theorem~\ref{theorem:tw-algo},
and so we do not prove it in detail here.  The basic idea is as follows.

We say that a boundary configuration $C$ is {\em viable} for a graph $G$ if there
exists some partial-3-manifold triangulation $\T$ with $\Gamma(\T) = G$
and with $C$ as its boundary configuration.
Our algorithm starts with an empty triangulation,
and then introduces tetrahedra and face identifications in a way that
essentially works from the leaves up to the root of the tree decomposition
of $G$.  For each subtree in the tree decomposition
we compute which configurations are viable for the
corresponding subgraph of $G$, and then propagate these configurations further up
the tree.  The running time at each node depends only on the number of
boundary faces, which is bounded in terms of the bag size and thereby
$\tw(G)$.

%\begin{prob} \textsc{valid face pairing graph}
%Given a connected 4-regular multigraph $G$, determine whether there exists a
%closed 3-manifold triangulation $\T$ such that $\Gamma(\T) = G$.
%\end{prob}

\subsection{A generalisation to simple properties} \label{sec:algorithm}

% While the earlier result in this section was the original aim of the research,
Here we generalise Theorem~\ref{theorem:tw-algo-plain}
to many other settings.
For this we define a \emph{simple property} of a partial 3-manifold
triangulation (see below).

We extend boundary configurations to include an extra piece
of data $\phi$ based on the partial triangulation that helps test our property.
For instance, if $p$ is the simple property that the triangulation
contains $\leq 3$ internal vertices,
then $\phi$ might encode the number of internal
vertices thus far in the partial 3-manifold triangulation
(here $\phi$ takes one of the values $0,1,2,3,\texttt{too\_many}$).

As before: for a simple property $p$,
we say that a boundary configuration $C$ is {\em $p$-viable}
for a graph $G$ if there exists some
partial-3-manifold triangulation $\T$ with property $p$, with
$\Gamma(\T)=G$ and with $C$ as its boundary configuration.

Shortly we solve the problem of testing whether a graph $G$
admits any closed 3-manifold triangulation with property $p$,
for any simple property $p$.  The basic idea is as before: for each
subtree of our tree decomposition of $G$, we compute all
viable
% (but not necessarily $p$-viable)
configurations
and propagate this information up the tree.

\begin{definition}[Simple property]
  \label{definition:simple}
  A boolean property $p$ of a partial-3-mani\-fold triangulation
  is called \emph{simple} if all of the following hold.
  Here all configurations have $\leq b$ boundary faces,
  and $f,g,h$ are some computable functions.
  \begin{enumerate}
    \item \label{definition:simple:values-p}
      The extra data $\phi$ in the boundary configuration
      satisfies $\phi \in P$ for some
      universal set $P$ with $|P| \leq f(b)$.
    \item We can determine whether a triangulation satisfies $p$ based only on
      its boundary configuration (including the extra data $\phi$).
    \item
      Given any viable configuration and a new face identification
      $\pi$ between two of its boundary faces, we can in $O(g(b))$ time
      test whether introducing this identification yields another
      viable configuration and, if so, calculate the corresponding
      value of $\phi$.
    \item Given viable configurations for two
      disjoint triangulations, we can in
      $O(h(b))$ time test whether the configuration for their union is
      also
      viable and, if so, calculate the corresponding value of $\phi$.
  \end{enumerate}
\end{definition}

\noindent
The four conditions above can be respectively interpreted as meaning:
\begin{enumerate}
  \item the upper bound on the number of viable configurations
    (including the data $\phi$) still depends on $b$
    but not the number of tetrahedra;
  \item we can still test property $p$ without examining the full triangulation;
  \item new face identifications can still be checked for $p$-viability in
    $O(g(b))$ time;
  \item configurations for disjoint triangulations can be combined
    in $O(h(b))$ time.
\end{enumerate}

\begin{example}\label{example:n-vertex}
Let $p$ be the property that a triangulation contains
at most $x$ internal vertices (i.e., vertices with links homeomorphic to a
2-sphere), for some fixed integer $x$.  Then $p$ is simple.

Here we define $\phi \in P = \{0,1,\ldots,x,\texttt{too\_many}\}$ to be the
number of vertices in our partial 3-manifold triangulation with
2-sphere links. This clearly satisfies conditions 1 and 2.
For condition 3: when identifying two faces together, a new vertex acquires a
2-sphere link if and only if the identification closes off all punctures
in the link (which we can test from the edge and vertex configurations).
Condition 4 is easily satisfied by summing $\phi$ over the disjoint
configurations.

The case when $x=1$ is highly relevant: much theoretical and
computational work has gone into
1-vertex triangulations of 3-manifolds
\cite{Jaco2002Algorithms,Matveev2007AlgorithmicTopology},
and these are of particular use when searching for 0-efficient
triangulations \cite{Jaco2003ZeroEfficient}.
\end{example}

%and Theorem \ref{theorem:tw-algo}
%shows that \vfpg{$p$} is fixed parameter tractable for $p$ being {\em the
%triangulation contains at most $x$ vertices} for any constant $x$.

We can now state the main result of this paper:

\begin{prob} \vfpg{$p$}{$G$}
Let $p$ be a simple property.  Given a connected 4-regular multigraph $G$,
determine whether there exists a closed 3-manifold triangulation $\T$
with property
$p$ such that $\Gamma(\T) = G$.
\end{prob}

\setcounter{c-save-main}{\arabic{theorem}}
\begin{theorem}\label{theorem:tw-algo}
Let $p$ be a simple property.
Given a connected 4-regular multigraph $G$ on $n$ nodes with treewidth $\leq
k$, and a corresponding tree decomposition with $O(n)$ nodes where each
bag has at most two children, we can solve
\vfpg{$p$}{$G$} in $O(n\cdot f(k))$ time for some computable function $f$.
\end{theorem}

Our requirement for such a tree decomposition is not restrictive:
Bodlaender \cite{Bodlaender1996} gives a fixed parameter tractable algorithm
to find a tree decomposition of width $\leq k$ for fixed $k$, and
Kloks \cite{Kloks1994Treewidth} gives an $O(n)$ time algorithm to
transform this into a tree decomposition where
each bag has at most two children.  The ``two children'' constraint can
be relaxed; we use it here to keep the proof simple.

%A full proof follows, we first outline the idea of the proof.
%%appears in Appendix~\ref{app:main}; the main ideas are
%%as follows.
%For each bag $\nu$ of the tree decomposition we define a
%corresponding subgraph $G_\nu$ of $G$, which contains precisely those nodes of
%$G$ that do \emph{not} appear in bags outside the subtree rooted at $\nu$.
%As before we use a dynamic programming approach, working from the
%leaves of the tree decomposition up to the root: for each $\nu$ we
%construct all viable
%% (but not necessarily $p$-viable)
%configurations for $G_\nu$,
%by combining the viable configurations at the
%child nodes of $\nu$ and analysing any new face identifications that
%might appear.
%We bound the running time at each $\nu$ by a
%function of the bag size, using the properties of
%Definition~\ref{definition:simple} and the observation that any
%partial triangulation admitted by
%$G_\nu$ must have $\leq 4(\tw(G)+1)$ boundary faces.
%
%Once we reach the root node of the tree decomposition, the
%final list contains a $p$-viable configuration if and only if
%$G$ admits a closed 3-manifold triangulation with property $p$.

\proof We will give this proof in three sections. First, we describe the
algorithm in detail. Then we show that the algorithm is correct. Lastly we
demonstrate the running time of the algorithm.

We begin, however, with some preliminaries.
Let the tree decomposition be $H$. Recall that in a tree decomposition,
each node of $H$ represents a collection of nodes of $G$. We will use
the term {\em bag} to
refer to a node of $H$, and {\em node} to refer to a node of $G$. Each node $w$
in $G$ will represent a corresponding tetrahedron~$\Delta_w$.

Arbitrarily choosing one bag of $H$ and make it the
root bag, so that the tree becomes a hierarchy of subtrees. For each bag $\nu$
in $H$, the subtree $H_\nu$ is defined as the subtree obtained by taking only
the bag $\nu$ and any bag which appears below $\nu$ in $H$.

We now define the subgraph $G_\nu$ of $G$, which contains precisely
those nodes of $G$ that appear \emph{only} in $H_\nu$.
In other words, node $w$ is in $G_\nu$ if and only if $w$
does \emph{not} appear in any bag in $H\setminus
H_\nu$.  The subgraph $G_\nu$ contains all corresponding arcs of $G$,
i.e., all arcs of $G$ that link nodes of $G_\nu$.

We first make the following observation.
If two children $\nu_i,\nu_j$ of some bag $\nu$ were to contain a common node
$w$, then since $H$ is a tree decomposition any such $w$ must also
appear in the bag $\nu$. Therefore
\emph{no two subgraphs $G_{\nu_i},G_{\nu_j}$ may
contain a common node representing a common tetrahedron}.

As a result, we can combine the boundary
configurations of children of $\nu$ by simply taking the union of the
configurations, as they correspond to disjoint triangulations.
The same process
can be used to extend some boundary configuration with the
boundary configuration of a new standalone tetrahedron.

\bigskip

{\noindent \underline{The algorithm:}}\quad
For each bag $\nu$, we will construct all boundary configurations
for $G_\nu$.  To achieve this, we take the following steps:

\begin{enumerate}
  \item \label{step:combine} Take every possible combination of configurations from the children of
    $\nu$, where each combination contains exactly one configuration from each
    child. We showed earlier that these must represent disjoint triangulations,
    and that for each combination we can construct the configuration of their union.
    %, so by Step 4 of Definition \ref{definition:simple} this can be done in
    %$O(h(k))$ time.
  \item \label{step:each-config} For each such combined configuration $C$:
    \begin{enumerate}
      \item \label{step:each-elt} For each element $w$ inside the bag $\nu$, if $w$ does not appear
    anywhere in a higher bag in $H$, then add the boundary configuration of a
    single standalone tetrahedron (corresponding to tetrahedron $\Delta_w$) to
    $C$. Again, the earlier observation shows that this is possible. Then:
    \begin{enumerate}
      \item \label{step:each-arc} For each arc $e$ in $G$ incident to the node $w$, if
        the other endpoint $w'$ of the arc $e$ is also in
        $G_\nu$, use Lemma \ref{lemma:add-one-id} to try to add each of the six
        possible corresponding face identifications to $C$ in turn
        (recall that these come from the six symmetries of a regular triangle). %, as per Step 3 of
        %Definition \ref{definition:simple}.
        For each
        viable (but not necessarily $p$-viable)
        configuration thus created, continue by recursing to Step
        \ref{step:each-elt} and taking
        the next element $w$.
    \end{enumerate}
  \item \label{step:found-config} If all elements of $\nu$ have been
  successfully processed in Step \ref{step:each-elt}
  then a viable
      configuration has been found. Store this as a viable boundary
      configuration for $G_\nu$.
\end{enumerate}
\end{enumerate}

If any bag contains no viable configurations, we immediately know that
there are no closed 3-manifold triangulations $\T$ satisfying \vfpg{$p$}{$G$}.

Once all configurations have been constructed, if the root bag
%If $\nu$ is the root bag (i.e. $G_\nu = G$) and $\nu$
contains any $p$-viable boundary configurations (by construction all
boundary configurations at the root node have empty edge and vertex
configurations), then there does exist some closed 3-manifold triangulation
$\T$ with property $p$ such that $\Gamma(\T) = G$. If, however, the root bag
contains no $p$-viable configuration then such a triangulation does not exist.
%Note
%the distinction being made between ``empty configurations'' being in a bag, and
%a bag not containing any configurations. An empty configuration may very well
%still be a $p$-viable configuration.

\bigskip

{\noindent \underline{Correctness:}}\quad
For each bag $\nu$, we have a corresponding graph $G_\nu$. If a closed
3-manifold triangulation $\T$ with property $p$ and $\Gamma(\T) = G$ exists,
then define $\T_\nu$ to be the partial-3-manifold triangulation constructed
by removing from $\T$ the tetrahedra
and face identifications which respectively represent nodes and arcs not
present in $G_\nu$. Each such $\T_\nu$ must be a partial-3-manifold
triangulation, and so each bag $\nu$ must have at least one viable
configuration.

If the root bag does contain some $p$-viable boundary configuration, then each
arc in $G$ has been through Step \ref{step:each-arc} in the algorithm and by Lemma
\ref{lemma:add-one-id} we know that each such configuration must represent a
partial 3-manifold triangulation with property $p$ (or possibly many such
triangulations). Since $G$ is 4-regular, we also know that these triangulations
can have no boundary faces, and so these partial 3-manifold triangulations must
in fact be closed 3-manifold triangulations with property $p$.

\bigskip

{\noindent \underline{Running time:}}\quad
We begin by showing that the number of configurations at each bag $\nu$ is
bounded by a function of $k$, but is independent of $n$.

Consider a boundary face $f$ of tetrahedron $\Delta$ in some triangulation
$\T_\nu$ represented by some configuration $C$ in $\nu$. In the graph $G$,
there must exist some arc $a$ which represents the identification of $f$ with
some other face $f'$ of some tetrahedron $\Delta'$. However, since $f$ is a
boundary face of $\T_\nu$, this must mean that the node representing $\Delta'$
must occur in some bag ``higher up'' in the tree decomposition; that is,
the node representing $\Delta'$ must occur in some
bag in $H\setminus H_\nu$. However, the nodes representing $\Delta$ and
$\Delta'$ must occur together in some bag (as they are the endpoints of arc
$a$), so by Definition \ref{definition:treedecomp} the node representing
$\Delta$ must occur in the bag $\nu$ itself. From this, it is easy to see that as $\nu$ has at
most $k+1$ elements, configurations at $\nu$ must represent triangulations with
at most $4(k+1)$ boundary faces. Therefore by Lemma
\ref{coro:num-bdry-configs} and Definition \ref{definition:simple}
the number of configurations is bounded by some computable function of $k$.

We now calculate the running time of each step in the algorithm.
The tree decomposition has $O(n)$ nodes, and at each node we go through all
three steps.
For
these steps we do not calculate exact upper bounds, as the functions tend to be
exponential or worse. We instead only show that the running time is bounded by
some function of $k$, which suffices to give the required result. For these
steps, $f$, $g$ and $h$ are used to denote some arbitrary computable functions
of $k$. They are not the same as those given in Definition
\ref{definition:simple}.

At Step \ref{step:combine}, we are combining configurations. The number of configurations in
each child bag is independent of $n$ by the above argument, and each bag has at
two children so for a single bag, Step \ref{step:combine} takes $O(f(k))$ time per bag.

Step \ref{step:each-config} takes each such configuration, and at Step
\ref{step:each-elt} extends the
configuration. Since each bag contains at most $k+1$ elements, Step
\ref{step:each-elt}  runs
at most $k+1$ times per bag. We know that $G$ is 4-regular, so by Lemma
\ref{lemma:add-one-id} and Definition \ref{definition:simple} there are at most
four distinct pairs of faces to identify,
% ^^ in what step? - bab
and thus Step \ref{step:each-arc} runs
in $O(g(k))$ time and therefore Step \ref{step:each-elt} likewise
runs in $O(k\cdot g(k))$ time. Step
\ref{step:found-config} is simply storing configurations. By condition
\ref{definition:simple:values-p} of Definition
\ref{definition:simple}, and the definition of boundary configurations, the
size of each configuration is independent of $n$.  Since the number of
configurations at each bag is also independent of $n$, Step
\ref{step:each-config} can therefore be
completed for each bag in $O(h(k))$ time.

Combining the above counts for each step for each of the $O(n)$ bags in the
tree decomposition gives the required result.
\qed

\section{Implementation and experimentation} \label{sec:impl}

The algorithm was implemented Java, using the treewidth library from
\cite{vanDijk2006}.
%Vertex,
%edge and boundary configurations were are all stored as hash maps, as although
Although our theoretical bound on the number of configurations
is extremely large (Lemma \ref{lemma:exact-num-bdry-configs}), we store all
configurations using hash maps to exploit situations
where in practice the number of viable configurations is
much smaller.  As seen below, we find that such a
discrepancy does indeed arise (and significantly so).

We also introduce another modification that yields significant speed
improvements in practice.
The algorithm builds up a complete list of
all viable configurations at each bag $\nu$ of the tree decomposition.
However, for an affirmative answer to the problem, only
a small subset of these may be required.
We take advantage of this as follows.

For any bag $\nu$ with no children, configurations are computed as normal.
Once a viable configuration is found, it is immediately propagated up
the tree in a depth-first manner.  This means that, rather than calculating
every possible viable configuration for every subgraph $G_\nu$,
the improved algorithm can identify a full triangulation
with property $p$ quickly and allow early termination.

We implemented the program with $p$ defined to be {\em one-vertex and possibly
minimal}, using criteria on the degrees of edges from \cite{Burton2004}. This
allowed us to compare both correctness and timing with the
existing software {\em Regina}
\cite{Regina}. We ran our algorithm on all 4-regular graphs on 4, 5 or 6 nodes
(see Table \ref{tab:regular})% in Appendix~\ref{app:expt} for details)
to verify correctness.
We see that the average
time to process a graph increases with treewidth, as expected. We also
see that the number of viable configurations is indeed significantly lower than
the upper bound of Lemma \ref{lemma:exact-num-bdry-configs}, as we had hoped.

%We finish with Table~\ref{tab:regular}, which contains experimental data
%for a thorough suite of smaller problems as discussed in Section~\ref{sec:impl}.

\begin{table}[h] \centering
\caption{Results from the original algorithm. From left to right, the columns
  denote the number of nodes in the graph, the treewidth of the
  graph, the number of distinct graphs with these parameters, the average
  running time of the algorithm on these graphs, and the largest number of
  configurations found at any bag, for any graph.}
\begin{tabular}{|c|c||c|c|c|}
\hline
$|(V(G))|$ & $\tw(G)$ & \# of graphs & Avg. run time (ms)& max$(|$configurations$|)$ \\
\hline \hline
4 & 1 & 1 & 680 & 2 \\
4 & 2 & 8 & 4036 & 7 \\
4 & 3 & 1 & 13280 & 17 \\ \hline
5 & 1 & 1 & 780 & 17 \\
5 & 2 & 22 & 13446 & 156 \\
5 & 3 & 4 & 29505 & 307 \\
5 & 4 & 1 & 94060 & 39 \\ \hline
6 & 1 & 1 & 890 & 17 \\
6 & 2 & 68 & 64650 & 1583 \\
6 & 3 & 25 & 346028 & 5471 \\
6 & 4 & 3 & 297183 & 1266 \\ \hline
\end{tabular}
\label{tab:regular}
\end{table}

{\em Regina} significantly outperforms our algorithm on all of these graphs,
though these are small problems for which asymptotic behaviour plays a
less important role.  What matters more is performance on
larger graphs, where existing software begins to break down.
%though we note that \emph{Regina} benefits from
%over 10 years of careful optimisation, compared to our initial
%proof-of-concept code.  What these experiments do indicate is that
%our algorithm is clearly viable, despite the enormous bounds of
%Lemma~\ref{lemma:exact-num-bdry-configs}.

We therefore ran a sample of 12-node graphs through our algorithm,
selected randomly from graphs which cause significant slowdown in
existing software.  This ``biased'' sampling was deliberate---our
aim is not for our algorithm to \emph{always} outperform existing
software, but instead to seek new ways of solving those difficult
cases that existing software cannot handle.
Here we do find success: our algorithm was
at times 600\% faster at identifying non-admissible graphs than {\em Regina},
though this improvement was not consistent across all trials.  More
detailed experiments will appear in the full version of this paper.

In summary: for larger problems,
our proof-of-concept code already exhibits far superior performance
for some cases that \emph{Regina} struggles with.
With more careful optimisation (e.g., for dealing with
combinatorial isomorphism), we believe
that this algorithm would be an important tool that complements
existing software for topological enumeration.

The full source code for the
implementation of this algorithm is available at
\url{http://www.github.com/WPettersson/AdmissibleFPG}.

\section{Applications and extensions}

We first note that our meta-theorem is useful: here we list several
simple properties $p$ that are important in practice,
with a brief motivation for each.
\begin{enumerate}
  \item \emph{One-vertex triangulations} are crucial for computation:
    they typically use very few tetrahedra,
    and have desirable combinatorial properties.
    This is especially evident with
    0-efficient triangulations \cite{Jaco2003ZeroEfficient}.
  \item Likewise, \emph{minimal triangulations} (which use the fewest
    possible tetrahedra) are important for both combinatorics and
    computation \cite{Burton2004,Burton2007}.
    Although minimality is not a simple property, it has many
    simple necessary conditions, which are
    used in practical enumeration software
    \cite{Burton2007,Matveev2007AlgorithmicTopology}.
  \item \emph{Ideal triangulation of hyperbolic manifolds} play a key role
    in 3-manifold topology.  An extension of Theorem~\ref{theorem:tw-algo}
    allows us to support several necessary conditions for hyperbolicity,
    which again are used in real software
    \cite{Callahan1999,Hildebrand1989}.
\end{enumerate}

%We now highlight another use of treewidth in topological algorithms.
%Although existing enumeration algorithms construct triangulations by
%identifying faces, there has been little research into determining a
%good ordering for which
%faces to identify first. By taking an ordering from our tree
%decomposition in Theorem \ref{theorem:tw-algo}, we have modified
%the software \emph{Regina}
%\cite{Regina} and improved the running time of a census
%%\footnote{Note that this census algorithm
%%is not fixed parameter tractable, as we only modified the ordering of the
%%tetrahedra.}
%on 9 tetrahedra from 274.47 CPU-minutes down to 65.906 CPU-minutes.
%%However, other experimental algorithms \cite{Pettersson2014edgelinkdecomp} show
%%much stronger running times again, and so
%Clearly more research into this would be beneficial.

Finally: a major limitation of all existing 3-manifold enumeration
algorithms is that they cannot ``piggyback'' on prior results for fewer
tetrahedra, a technique that has been remarkably successful in other
areas such as graph enumeration \cite{mckay98-generation}.
This is not a simple oversight: it is well known that we cannot build all
``larger'' 3-manifold triangulations from smaller 3-manifold triangulations.
The techniques presented
here, however, may allow us to overcome this issue---we can modify the algorithm
of Theorem \ref{theorem:tw-algo} to store entire families of triangulations at
each bag of the tree decomposition. We would lose fixed parameter
tractability, but for the first time we would be able to cache and reuse
partial results across different graphs and even different numbers of
tetrahedra, offering a real potential to extend census data well beyond
its current limitations.
%algorithm is no longer fixed parameter tractable, but it is easy to see that it
%is possible for two distinct graphs $G$ and $G'$ to have, in the notation of the
%algorithm, similar graphs $G_\nu$ and $G'_{\nu'}$. If $G_\nu$ and $G'_{\nu'}$
%are indeed isomorphic, then the results obtained in when running the algorithm
%on $G$ can be used when running the algorithm on $G'$, and vice-versa.

\bibliography{mybib}
\bibliographystyle{amsplain}

%\newpage
%\appendix
%\section*{Appendix}

%\begin{theorem}[\ref{theorem:cutwidth-algo}]
%Let $G$ be a connected 4-regular multigraph on $n$ nodes with cutwidth $\leq
%k$, and for which a corresponding layout of nodes is known, and let $p$ be a
%simple property. Then \vfpg{$p$} can be solved for $G$ in $O(n\cdot f(k))$ time.
%\end{theorem}
%
%\proof 
%\qed
%
%\begin{theorem}[\ref{theorem:tw-algo}]
%Let $G$ be a connected 4-regular multigraph on $n$ nodes with treewidth $\leq
%k$, and for which a corresponding tree decomposition with $O(n)$ nodes is
%known. Let $p$ be a simple property. Then \vfpg{$p$} can be solved for $G$ in
%$O(n\cdot f(k))$ time.
%\end{theorem}
%
%\proof
%\qed

%\section{Additional examples} \label{app:examples}
%
%Here we give an additional example of edge and vertex configurations,
%as described in Section~\ref{sec:config}.
%
%
%\section{Proof of Lemma~\ref{lemma:add-one-id}}
%\label{app:add-one-id}
%
%
%\section{Proof of Theorem~\ref{theorem:tw-algo}} \label{app:main}
%
%\setcounter{theorem}{\arabic{c-save-main}}

%\section{Experimental data} \label{app:expt}

\end{document}